\documentclass[11pt]{article}

\input{amssym.def}
\input{amssym}
\usepackage{eufrak}

\newtheorem{theorem}{Theorem}

\newtheorem{definition}[theorem]{Definition}
\newtheorem{example}[theorem]{Example}
\newtheorem{lemma}[theorem]{Lemma}
\newtheorem{proposition}[theorem]{Proposition}
\newtheorem{remark}[theorem]{Remark}

\def\pa{{\partial}}
\def\hl{\hat{l}}
\def\hL{\hat{L}}
\def\DD{\cal{D}}

\def\un{\underline}
\def\qq{q^{-1}}
\def\PP{{\cal P}}

\def\Tr{\mathrm{Tr}}

\def\hLL{\hat{{\cal L}}}

\def\MM{{\cal{M}}}

\def\RR{R^{-1}}
\def\DD{{\cal{D}}}

\def\JJ{J^{-1}}

\def\de{\delta}

\def\ot{\otimes}
\def\C{{\Bbb C}}

\def\vv{V^{\otimes 2}}

\def\ov{\overline}

\def\la{{\lambda}}

\def\be{\begin{equation}}
\def\ee{\end{equation}}

\begin{document}

\makeatletter
\renewcommand{\theequation}{{\thesection}.{\arabic{equation}}}
\@addtoreset{equation}{section} \makeatother

\title{Wick theorem and matrix Capelli identities for quantum differential operators on Reflection Equation Algebras}
\author{\rule{0pt}{7mm} Dimitry
Gurevich\thanks{gurevich@ihes.fr}\\
{\small\it Higher School of Modern Mathematics,}\\
{\small\it Moscow Institute of Physics and Technology}\\
\rule{0pt}{7mm} Pavel Saponov\thanks{Pavel.Saponov@ihep.ru}\\
{\small\it
National Research University Higher School of Economics,}\\
{\small\it 20 Myasnitskaya Ulitsa, Moscow 101000, Russian Federation}\\
{\small \it and}\\
{\small \it
Institute for High Energy Physics, NRC "Kurchatov Institute"}\\
{\small \it Protvino 142281, Russian Federation}\\
\rule{0pt}{7mm} Mikhail Zaitsev\thanks{mrzaytsev@edu.hse.ru}\\
{\small\it National Research University Higher School of Economics,}\\
{\small\it 20 Myasnitskaya Ulitsa, Moscow 101000, Russian Federation}}

\maketitle

\begin{abstract}
Quantum differential operators on Reflection  Equation Algebras, corresponding to Hecke symmetries $R$ were introduced in previous publications. In the present paper we are
mainly interested in quantum analogs of the Laplace and Casimir operators, which are invariant with respect to the action of the Quantum Groups $U_q(sl(N))$, provided $R$ is
the Drinfeld-Jimbo $R$-matrix. We prove that any such an operator maps the central characteristic subalgebra of a Reflection Equation
algebra into itself. Also, we define the notion of normal ordering for the  quantum differential operators and prove an analog of the Wick
theorem for the product of partially ordered operators. As an important corollary we find a set of universal matrix Capelli identities
generalizing the results of \cite{Ok2} and  \cite{JLM}. Besides, we prove that the normal ordered form of any central
differential operator from the characteristic subalgebra is also a central differential operator.
\end{abstract}

{\bf AMS Mathematics Subject Classification, 2020: 17B37, 81R50}


\section{Introduction}
In this paper we introduce and study quantum analogs of  differential operators defined on a class of quantum matrix algebras called Reflection Equation (RE) algebras. We use
the term {\em quantum analog} in a wide sense of the word, not restricting ourselves to objects, related to the quantum groups
$U_q(sl(N))$. The latter objects ($R$-matrices, differential operators etc) are called {\em standard} in what follows.

The starting point for constructing an RE algebra is a {\it Hecke symmetry} $R$. Any such a symmetry is a linear operator $R\in \mathrm{End}(V^{\otimes 2})$, $V$ being a complex
vector space $\dim_{\,\Bbb C}V = N$, which is subject to a {\it braid relation} in the algebra $\mathrm{End}(V^{\otimes 3})$
\be
(R\ot I)(I\ot R)(R\ot I)=(I\ot R)(R\ot I)(I\ot R),
\label{braid}
\ee
and, in addition, obeys the {\it Hecke condition}
\be
(q\,I\otimes I-R)(q^{-1}I\otimes I+R)=0,\qquad q\in \C\setminus \{\pm 1,0\}.
\label{Hecke}
\ee
Hereafter, $I$ stands for the identity operator on $V$ or the $N\times N$ unit matrix. The numeric parameter $q$ is assumed to be {\it generic}, i.e. $q^{n}\not=1$ for all positive integers $n\in \Bbb{Z}_+$.

\begin{remark}\rm
If in (\ref{Hecke}) one sets $q=1$, then the corresponding operator $R$ is called an {\it involutive symmetry}.
The well known example of an involutive symmetry is the usual flip $P=\mathrm{End}(V^{\otimes 2}):$
$$
P\triangleright (u\otimes v) = v\otimes u,\qquad \forall\,u,v\in V.
$$
From now on the symbol $\triangleright$ denotes the action of a linear operator.

Below, we suppose that any Hecke symmetry in question allows a limit $q \to 1$ and turns into  an involutive symmetry in this limit.
If, in addition, $\lim_{q\to 1}R = P$, then $R$ will be referred to as a deformation of the flip $P$. The most important example
of such a Hecke symmetry is the Drinfeld-Jimbo $R$-matrix connected with the quantum group $U_q(sl(N))$.
\end{remark}

Given a Hecke symmetry $R$, the corresponding RE algebra is defined as a unital associative algebra $\MM(R)$ finitely generated
by the set of entries of an $N\times N$ matrix $M=\|m_i^j\|$, which satisfies the relation:
\be
R\,(M\ot I)R\,(M\ot I)-(M\ot I)R\, (M\ot I)R=0.
\label{RE}
\ee
Note that if $R=P$, then $\MM(P) = \mathrm{Sym}(gl(N))$.

The RE algebras possess a lot of remarkable properties. In particular, they admit introducing  of analogs of the partial derivatives in the
generators $m_i^j$ and analogs of the Weyl-Heisenberg
algebras. The latter analogs can be introduced in terms of the so-called {\it quantum doubles}. By a quantum double $(A,B)$ we mean
a couple of associative algebras $A$ and $B$ equipped with a {\em permutation map} $\sigma: A\ot B\to B\ot A$, compatible with their
algebraic structures (see \cite{GS1} for more details and precise definitions). The map $\sigma$ plays the role of the Leibniz
rule for elements of the algebra $A$, treated as vector fields (in particular,  quantum partial derivatives), which act on  elements of the
algebra $B$. Note that the classical Leibniz rule is universal and is valid for all vector fields, whereas in the quantum case the permutation
map essentially depends on given algebras $A$ and
$B$. Observe that by contrast with the approach of the paper \cite{IP}, in our setting the quantum double consists of two RE algebras in
different realizations\footnote{Considered in \cite{IP} was the double of an RE algebra and an algebra of quantized functions on a group
(the so-called RTT algebra).}.

In the quantum double, we are dealing with, the role of a function subalgebra $B$ is played by the RE algebra $\MM(R)$. The algebra $A$ is
realized as the RE algebra $\DD(R^{-1})$, generated by entries of the  matrix $D=\|\pa_i^j\|$ which is subject to the relation:
\be
R^{-1}(D\ot I)R^{-1}(D\ot I)-(D\ot I)R^{-1} (D\ot I)R^{-1}=0.
\label{RE-D}
\ee
The permutation relations among generators of the algebras $\DD(\RR)$ and $\MM(R)$ are exhibited in (\ref{raz}). This
relations allow us to define a linear action of the subalgebra $\DD(R^{-1})$ on subalgebra $\MM(R)$ which is compatible both with
the algebraic structure of  $\DD(R^{-1})$ and that of $\MM(R)$. Due to this fact, the generators $\partial_i^j$ of the RE algebra
$\DD(R^{-1})$ are treated as {\it quantum derivatives}. Note that in the standard case the corresponding
double $(\DD(R^{-1}),\MM(R))$  is a quantization (deformation) of the  $gl(N)$-type Weyl-Heisenberg algebra of differential
operators in partial derivatives $\partial/\partial m_i^j$ with polynomial coefficients in coordinates $m_i^j$.

An important property of the quantum double $(\DD(R^{-1}),\MM(R))$  is that the entries of the matrix $\hL=MD = \|\hat l_i^j\|$
generate the so-called modified RE algebra $\hLL(R)$. The permutation relations among its generators are determined by
the following matrix equality:
\be
R\, (\hL\ot I)R\, (\hL\ot I)-(\hL\ot I)R\, (\hL\ot I) R=R\, (\hL\ot I)-(\hL\ot I)R.
\label{mRE}
\ee
If $R=P$, the algebra $\hLL(P)$ is isomorphic to the universal enveloping algebra $U(gl(N))$. In this case the generators
$\hl_i^j=\sum_k m_i^k \pa^j_k$ of the  algebra $U(gl(N))$
are represented by the Euler type vector fields, acting on the commutative  algebra
$\MM(P) = \mathrm{Sym}(gl(N))\simeq {\Bbb C}[m_i^j]$.
In general,  we get  quantum analogs of these vector fields acting on the noncommutative RE algebra
$\MM(R)$, whereas  a new quantum Leibniz rule for them is given by the corresponding permutation relations.

In the classical case (i.e. while $R=P$)  the elements $\Tr D^{k}$ and $\Tr \hL^k$ are invariant under the adjoint action of the group $GL(N)$
and are referred to as Laplace and Casimir operators respectively. Also, the normally ordered operators
$$
W^{(\mathbf{k})} =\, :\!\Tr(\hL^{k_1})\dots \Tr (\hL^{k_p})\!: \qquad \mathbf{k}=(k_1,\dots,k_p), \quad k_i\in {\Bbb Z }_{\ge 0}
$$
are of interest. Here, the colons stand for the usual normal ordering
\begin{equation}
:\!\partial_k^s m_i^j\!:\, = m_i^j \partial_k^s.
\label{cl-ord}
\end{equation}
Being properly normalized, the operators $W^{(\mathbf{k})}$ are
called cut-and-join ones. They play an important role in combinatorics and integrable system theory.

Note that numerous aspects of classical theory can be extended to the RE algebras and their doubles --- algebras of quantum
differential operators. Below, we define the quantum analogs of the Laplace and Casimir operators in a way, similar to the classical
patterns.  The peculiarity of the quantum case is that the usual matrix trace $\Tr$, entering the
definitions of invariant elements, is replaced by its quantum analog  $\Tr_R$, which is well defined for any {\em skew-invertible}
Hecke symmetry $R$. In the standard case the quantum Laplace and Casimir operators are invariant with respect to the action of
the quantum group $U_q(sl(N))$.

We define the notion of normal ordering of quantum derivatives $\partial_i^j$ and ``coordinates" $m_i^j$ which is a deformation of
the classical rule (\ref{cl-ord}). Emphasize, that the relation $\hL=MD$ enables us to define the ordered form of Casimir operators in the
quantum case. We prove that any normally ordered quantum differential operator from the central characteristic subalgebra
of the algebra $\hLL(R)$ is a central differential operator. In this way, the quantum analogs of the cut-and-join operators
$W^{(\mathbf{k})}$ can be easily defined.

Moreover, we generalize on the quantum double $(\DD(R^{-1}),\MM(R))$ the Wick formula,
relating differential operators and their normally ordered forms. This formula entails a set of
new matrix Capelli identities, generalizing those from \cite{Ok2} and \cite{JLM}. Thus, the main ingredients of the paper
\cite{Ok2} --- the symmetric groups and the enveloping algebras $U(gl(N))$ are replaced by their quantum analogs
--- the Hecke and RE algebras. As for the matrix Capelli identities, in the recent paper \cite{JLM} they were proved in a
particular case of the standard Hecke symmetry. By contrast with \cite{JLM}, our version of the Capelli identities is universal
in the sence that they are not attached to a concrete Young diagram and are valid for {\it any} skew-invertible Hecke symmetry
$R$ including the supersymmetric $GL(m|n)$ type generalizations.  Also, note that quantum analogs of immanants from \cite{Ok1}
can be defined not only in the standard case as was done in \cite{JLM} but for any RE algebra (see \cite{Z}).

The paper is organized as follows. In the next section we remind some basic constructions related to the Hecke symmetries and the
corresponding RE algebras. In Section \ref{sec:3} we introduce the quantum derivatives and define generalized quantum
Laplacians on these algebras. We prove that the action of the generalized Laplace operators maps the central
characteristic subalgebra into itself. In Section \ref{sec:4} we introduce the generalized Casimir operators by means of the ``quantum
Euler-type vector fields". In Section \ref{sec:5} we define the normal ordering for quantum differential operators and prove a version
of the Wick theorem which allows to transform the partially ordered product of quantum differential operators into totally 
ordered form. Section \ref{sec:6} is devoted to derivation of the set of universal matrix Capelli identities, as well as to the proof of
centrality of the normally ordered Casimir operators.

\section{Preliminaries}
\label{sec:2}

In this section we introduce some basic notation and give definitions of objects and constructions used in what follows.

Recall that for any given integer $k\ge 2$ the {\it Hecke algebra} $H_k(q)$ is a unital associative algebra over the complex field
${\Bbb C}$ generated by Artin's generators $\tau_1,\dots, \tau_{k-1}$, which are subject to the set of relations:
$$
\begin{array}{lcl}
	\tau_i\,\tau_{i+1}\, \tau_i=\tau_{i+1}\,\tau_i\,\tau_{i+1},&\quad&1\le i\le k-2\\
	\rule{0pt}{5mm}
	\tau_i\, \tau_j=\tau_j\,\tau_i,&&|i-j|\geq 2\\
	\rule{0pt}{5mm}
	(q\, 1-\tau_i)(\qq\, 1+\tau_i)=0,&&q\in \C\setminus \{\pm 1, 0\},
\end{array}
$$
where $1$ stands for the unit element of the algebra. As is known, the Hecke algebra $H_k(q)$ is finite dimensional and for a
generic $q$ is isomorphic to the group algebra ${\Bbb C}[S_k]$ of the symmetric group $S_k$. It is convenient to define $H_1(q)$
as an algebra generated by the only unit element. Thus, we have $H_1(q)\simeq {\Bbb C}$.

Below we need the so-called {\it Jucys-Murphy elements} $j_r$, $1\le r\le k$, which generate a maximal commutative subalgebra in
$H_k(q)$:
\be
j_1 = 1,\quad j_r = \tau_{r-1}\,j_{r-1}\,\tau_{r-1}, \quad 2\le r\le k.
\label{JM}
\ee

Observe that any Hecke symmetry $R$ defines the {\it $R$-matrix representation} $\rho_R$ of the Hecke algebra $H_k(q)$ in
the space $V^{\ot k}$:
$$
\rho_R: H_k(q)\rightarrow \mathrm{End}(V^{\otimes k}),\qquad \forall \,k\ge2.
$$
The representation $\rho_R$ is completely defined by the images $\rho_R(\tau_i)\in \mathrm{End}(V^{\otimes k})$
of the generators $\tau_i$:
\be
\rho_R(\tau_i) = R_i :=R_{i\, i+1}= I^{\ot (i-1)}\ot R\ot I^{\ot (k-i-1)},\qquad 1\le i \le k-1.
\label{rho-R}
\ee
Also, we denote $J_r:=\rho_R(j_r)$.

From now on, we assume any Hecke symmetry $R$ to be {\it skew-invertible}. This means that there exists an operator $\Psi:\vv\to \vv$
such that
\be
\Tr_2 R_{12}\, \Psi_{23} = \Tr_2\Psi_{12}R_{23} = P_{13}.
\label{psi}
\ee
Recall that $P$ stands for the usual flip.

Let us fix a basis $\{x_i\}_{1\le i\le N}$ in the space $V$ and the corresponding bases  $\{x_{i_1}\ot\dots  \ot x_{i_k}\}$ in the spaces
$V^{\ot k}$, $k\ge 2$. Then the operator $\Psi$ can
be identified with its $N^2\times N^2$ matrix $\|\Psi_{ij}^{\,rs}\|$ in the above tensor basis of the space $V^{\otimes 2}$. Consider
the $N\times N$ matrrix $C$ with the following matrix elements:
$$
C_i^j=\sum_{k=1}^N\Psi_{ik}^{jk}.
$$
The matrix $C$ allows one to define the {\it $R$-trace} of any $N\times N$ matrix $X$:
\be
\Tr_R X =\Tr (C X).
\label{R-tr}
\ee
An immediate consequence of (\ref{psi}) is the following important property:
\begin{equation}
\Tr_{R(2)} R_{12} = \Tr_2(C_2R_{12}) = I_1.
\label{trR}
\end{equation}

Now, we describe a way of constructing some central elements of the RE algebra $\MM(R)$ defined in (\ref{RE}), which belong to
its {\it characteristic subalgebra}. Also, we recall the definition of the quantum Schur functions (polynomials) and power sums as 
particular examples of these elements. 

For this purpose it is convenient to introduce the following
matrix ``copies'' of the generating matrix $M$:
\be
M_{\ov 1}=M_1 = M\otimes I^{\otimes (k-1)},\qquad M_{\ov {r+1}}= R_{\, r} M_{\,\ov r}\,R_{\, r}^{-1},\qquad r\ge 1.
\label{M-cop}
\ee
In virtue of definition (\ref{rho-R}) the matrices $M_{\,\ov{r}}$ for $1\le r\le k$ represent different embeddings of the quantum
generating matrix $M$ into the space $(\mathrm{Mat}_N)^{\otimes k}$ of $N^k\times N^k$ matrices. Below, we do not fix the
concrete value of the integer $k$ just assuming $k$ to be sufficiently large so that all the matrix formulae make sense.

In analogy with (\ref{rho-R}) we will often consider embeddings of arbitrary $N\times N$ matrices into the space 
$(\mathrm{Mat}_N)^{\otimes k}$ and will use a similar notation:
$$
X_i = I^{\otimes (i-1)}\otimes X\otimes I^{\otimes (k-i)},\quad \forall\, X\in\mathrm{Mat}_N.
$$
Note that for $X=I$ all these embeddings coincide: $I_i = I^{\otimes k}$ for any $1\le i\le k$.

The following important theorem was proved in  \cite{IOP}.

\begin{theorem} {\bf \cite{IOP}}
\label{prop:2}
Let $z\in H_n(q)$ be an arbitrary element. Then the $n$-th order homogeneous polynomial in generators $m_i^j$ of the RE algebra
$\MM(R)$
\rm
\be
ch_n(z)=\Tr_{R(1\dots n)}\left(\rho_R(z) M_{\ov{1\to n}} \,\right) = \Tr_{R(1\dots n)} \left(M_{\ov{1\to n}}\,\rho_R(z)\right),\quad
\forall\,n\ge 1
\label{sem}
\ee
\it
is central in the RE algebra $\MM(R)$. Here $ M_{\ov{1\to n}}:= M_{\ov{1}}M_{\ov 2}\dots M_{\ov{ n}}$, and
$$
\Tr_{R(1\dots n)}(X) = \Tr_1(\Tr_2(\dots(\Tr_n(C_1C_2\dots C_n X))\dots )),\qquad \forall\,X\in(\mathrm{Mat}_N)^{\otimes n}.
$$
\end{theorem}

In \cite{IOP} the map
$$
ch_n: H_n(q)\to Z(\MM(R)),\quad z\mapsto ch_n(z),
$$
was called {\em characteristic }. Here the notation $Z(A)$ stands for the center of the algebra $A$. Moreover, the direct
sum of the images of all maps $ch_n$ for  $n\in {\Bbb Z}_+$ is a central subalgebra in $\MM(R)$ referred to as the {\it characteristic
 subalgebra} of the RE algebra $\MM(R)$ (see \cite{IOP} for detail).

\begin{remark}\rm
For a general skew-invertible Hecke symmetry $R$ the characteristics subalgebra is at least a subset of the center of RE algebra
$\MM(R)$.  When investigating the center of the quantum group $U_q(sl(N))$, the authors of \cite{FRT} constructed an embedding
of the standard RE algebra $\MM(R)$ into $U_q(sl(N))$. This embedding was used in \cite{JLM} for proving the Capelli identities in
the particular case related to the standard Hecke symmetry.
\end{remark}

Consider now some specific elements of the characteristic subalgebra which we need in what follows.  For a given integer $k$ we chose the so-called Coxeter element
$\tau_{k-1}\tau_{k-2}\dots \tau_1\in H_k(q)$. Its image under the characteristic map $ch_k$ reads:
\be
p_k(M):=\Tr_{R(1\dots k)}\left(\rho_R(\tau_{k-1}\tau_{k-2}\dots \tau_1) M_{\ov{1\to k}}\,\right)=\Tr_{R(1\dots k)}\left( R_{k-1} R_{k-2}\dots R_1\,M_{\ov{1\to k}}\,\right).
\label{pws}
\ee
The polynomials $p_k(M)$, $k\ge 1$ are called the quantum {\em power sums}. Note that the power sums (\ref{pws}) can be reduced  to $p_k(M)=\Tr_R M^k$. In this form
they are similar to the classical ones, corresponding to the case $R=P$. The only difference is that the $R$-trace is used instead of the usual one.

For any partition $\la=(\la_1\geq \dots \geq\la_s)$, where $\la_i$ are non-negative integers, we introduce the corresponding symmetric function (polynomial) $p_\lambda(M)$
\be
p_{\la}(M) = p_{\la_1}(M)\dots p_{\la_s}(M).
\label{gen-ps}
\ee
Since all factors in this product are central in the algebra $\MM(R)$, their order does not matter. We  call the symmetric function $p_\lambda(M)$ the power sum, corresponding
to the partition $\la$. Thus, the power sums $p_k(M)$, $k\ge 1$ correspond to one-row partitions $\la = (k)$.

Another important set of central elements is formed by the {\it quantum Schur polynomials} $s_\lambda(M)$, associated with partitions $\lambda\vdash n$, $n\ge 1$. The polynomials
$s_\lambda$ were defined in \cite{GPS1} for more general class of quantum matrix algebras than the RE algebras. In the particular case
 of the RE algebra $\MM(R)$ the quantum
Schur polynomial $s_\lambda(M)$ is defined by the formula:
\be
s_\la(M)=\Tr_{R(1\dots n)}\left( \rho_R(e^\la_{T})M_{\ov{1\to n}}\,\right),\qquad \lambda\vdash n.
\label{schu}
\ee
Here $T$ is one of the standard Young tables corresponding to the Young diagram of the partition $\lambda$, while $e^\lambda_T\in H_n(q)$ is a primitive idempotent of the Hecke
algebra (see, for example, the review \cite{OP} for more detail). As was shown  in \cite{GPS1}, the right hand side of (\ref{schu}) depends
 only on the diagram $\lambda$ and does not depend on the table $T$.
The polynomials (\ref{schu}) satisfy the Littlewood-Richardson rule
$$
s_\lambda(M)s_\mu(M) = \sum_{\nu}C_{\lambda\mu}^\nu s_\nu(M),
$$
with classical coefficients $C_{\lambda\mu}^\nu$ \cite{GPS1}. This property justifies the term ``Schur polynomials'' for the elements (\ref{schu}).

\begin{remark} \rm
The connection of quantum polynomials $p_\lambda(M)$ and $s_\lambda(M)$ with power sums and Schur functions of the classical theory of symmetric functions becomes more clear
after a parameterization of the quantum polynomials by the ``eigenvalues" of the generating matrix $M$. This parametrization was introduced and studied in \cite{GPS1, GPS2, GPS4}. In \cite{GS2} it was treated to be a quantum analog of the Harish-Chandra map.
\end{remark}

\section{Quantum derivatives and generalized Laplace operators}
\label{sec:3}

Let us briefly remind the construction of a quantum double from \cite{GS1} since it plays the central role in the subsequent considerations.

By a quantum double we mean a couple $(A,B)$ of two  associative algebras $A$ and $B$, equipped with a {\it permutation map}
$$
\sigma: A\ot B\to B\ot A,
$$
which preserves the algebraic structures of $A$ and $B$. If these algebras are defined via relations on their generators, this property means that the ideals, generated by the
relations, are preserved by the map $\sigma$ (see \cite{GS1} for detail).

Often it is more convenient to use the so-called {\em permutation relations}
\be
a\ot b=\sigma(a\ot b), \quad  a\in A,\,\, b\in B. \label{pere}
\ee

If the algebra $A$ admits a one-dimensional representation $\varepsilon:A\to \C$, then it is possible to define an action $\triangleright : A\otimes B\rightarrow B$ of the algebra
$A$ on $B$ by setting:
\be
a\triangleright b:=(\mathrm{Id}\ot \varepsilon)\sigma (a\ot b), \quad \forall\, a\in A, \,\, b\in B.
\label{AB-act}
\ee
Here, as usual we identify the elements $b$ and $b\ot 1$. Due to the properties of $\sigma$, such an action will be a representation of $A$ in the algebra $B$ (see \cite{GS1}).

Our basic example of the quantum double is provided by  the RE algebras $A=\DD(\RR)$ and $B=\MM(R)$ defined by the relations (\ref{RE-D}) and (\ref{RE}).
Let us  introduce the permutation relations as follows \cite{GPS5}:
\be
D_1R_{12}M_1R_{12}=R_{12}M_1R_{12}^{-1}D_1+ R_{12}1_B1_A.
\label{raz}
\ee

The entries of the matrix $D=\|\pa_i^j\|_{1\leq i,j\leq N}$ acquire the meaning of operators if we apply the trivial one-dimensional representation of the RE algebra  $A = \DD(\RR)$:
$$
\varepsilon(1_A)=1,\qquad  \varepsilon(\pa_i^j)=0,\quad 1\le i, j\le N.
$$
Using the general formula (\ref{AB-act}) with this map $\varepsilon$ we get:
\be
D_1\triangleright R_{12}M_1R_{12} = R_{12}M_1R_{12}^{-1}\varepsilon(D_1)+ R_{12}1_B\varepsilon(1_A) = R_{12} 1_B.
\label{D-M}
\ee
Below we shall omit the symbols of the unit elements $1_A$ and $1_B$.

\begin{remark}\rm
Note, that the above formula for action of $D$ contains the summation over matrix indices, so it actually describes the actions of some linear combinations of $\partial_i^j$ on some
other linear combinations of $m_r^s$. Nevertheless, due to the invertibility and skew-invertibility of $R$, the action (\ref{D-M}) can be transformed into the following equivalent form:
$$
D_1\triangleright M_2 = \Tr_{0}(\Psi_{02})P_{12} \quad\Leftrightarrow \quad \partial_i^j\triangleright m_k^s = \delta_i^s\,B_k^j,\quad B_i^j = \sum_{k=1}^N\Psi_{ki}^{\,kj}.
$$
So, in fact, formula (\ref{D-M}) allows one to find the action of any given $\partial_i^j$ on any given $m_k^s$.
\end{remark}

Note that if  $R=P$, the algebras $\MM(P)$ and $\DD(P^{-1})$ become commutative and the action (\ref{D-M}) takes the form:
$$
D_1\triangleright P_{12}M_{1}P_{12} = D_1\triangleright M_2 = P_{12}\quad \Leftrightarrow \quad \pa_i^j\triangleright m_k^s=\de_i^s\, \de_k^j.
$$
This motivates us to consider the generators of the subalgebra $A=\DD(\RR)$ as quantum analogs of the usual partial derivatives 
$\partial_i^j = \partial/\partial m_j^i$. Also, we treat
the quantum double $(\DD(\RR), \MM(R))$ as a quantum counterpart of the usual Weyl-Heisenberg algebra of $gl(N)$ type.

Introduce the following elements which are central in the subalgebra $\DD(R^{-1})$ of the quantum double (\cite{IOP}, see also 
Proposition \ref{prop:2}): 
\begin{equation}
D_Q^{(m)}=\Tr_{R(1\dots m)}\left( Q(R_1,\dots, R_{m-1})\, D_{\,\ov{1\rightarrow m}}\right),\qquad m\ge 1,
\label{gen-Lap}
\end{equation}
where $Q(R_1,\dots,R_{m-1})$ is an arbitrary polynomial in $R_i$, $1\leq i\leq m-1$. Taking $Q = R_{m-1}\dots R_1$ we
get a set of power sums:
\begin{equation}
p_m(D) =  \Tr_{R(1\dots m)}\left(R_{m-1}\dots R_1 D_{\,\ov{1\rightarrow m}}\right)  = \Tr_R D^{\,m},\qquad  m\ge 1,
\label{Lap}
\end{equation}
where  $D_{\,\ov{1\rightarrow m}} = D_{\,\ov 1}D_{\,\ov 2}\dots D_{\,\ov m}$.

\begin{remark}\rm
Strictly speaking, we should use the $\Tr_{R^{-1}}$ operation in the RE algebra
$\DD(R^{-1})$. But  if $R^{-1}$ is also skew-invertible (we always assume this to be the case),  then $C_{R^{-1}}$ differs from $C_R$
by a scalar nonzero multiplier \cite{Og}. So, when constructing the central elements in the RE algebra $\DD(R^{-1})$, we can use the 
matrix $C_R$ and the $\Tr_R$ operation instead of $\Tr_{R^{-1}}$.
\end{remark}

\begin{definition} \rm 
The quantum differential operators, corresponding to elements (\ref{Lap}) and (\ref{gen-Lap})
with respect to the action (\ref{D-M}), are called respectively the {\em quantum Laplace} and the {\it generalized quantum Laplace} 
operators (Laplacians) on the RE subalgebra $\MM(R)$ of the double $(\DD(R^{-1}), \MM(R))$.
\end{definition}

The following theorem holds true.

\begin{theorem}\label{theor:8}
The action of any generalized quantum Laplace operator $D_Q^{(m)}$ maps the characteristic subalgebra of the
RE algebra $\MM(R)$ into itself.
\end{theorem}

\smallskip

\noindent
{\bf Proof.} To prove the claim we show that the action of any generalized Laplace operator $D_Q^{(m)}$ on an arbitrary element
of the characteristic subalgebra results in an element of the characteristic subalgebra. To do so we need the formula for the action
of an arbitrary $m$-th order monomial in quantum derivatives $D$ on an arbitrary $k$-th order monomial
in $M$, that is we have to find the result of the action:
$$
D_{\,\ov{1\rightarrow m}}\triangleright M_{\,\ov{m+1}}M_{\,\ov{m+2}}\dots M_{\,\ov{m+k}}
$$
for any given pair of positive integers $m\le k$.

It is convenient to introduce a shorthand notation:
$$
R_{i\to j}^{\pm 1} = \left\{
\begin{array}{lcl}
R_i^{\pm 1}\,  R_{i+1}^{\pm 1} \dots  R_j^{\pm 1}&\quad & j\geq i\\
\rule{0pt}{6mm} R_i^{\pm 1}\,  R_{i-1}^{\pm 1}\dots  R_j^{\pm 1}&& j<i\\
\end{array}\right. .
$$
Upon multiplying the relation (\ref{raz}) by $R_1^{-2}$ from the right, we get
\be
D_1M_{\ov 2}=M_{\ov 2}D_1 R^{-2}_1+R_1^{-1}.
\label{rela}
\ee
Then, taking into account the definition (\ref{M-cop}) of $M_{\,\ov r}$ one generalizes (\ref{rela}) to the form:
$$
D_1M_{\,\ov r}=M_{\,\ov r} D_1 R_{r-1\to 2}\, R_1^{-2} \RR_{2\to r-1}+ R_{r-1\to 2}\,\RR_{1\to r-1}.
$$
With a simple recursion this allows us to come to the permutation relations of $D$ with an arbitrary monomial in generators $M$:
\begin{eqnarray*}
D_1 M_{\,\ov 2}\dots M_{\,\ov {k+1}}=M_{\,\ov 2}\dots M_{\,\ov{k+1}} \,D_1 &&\hspace*{-6.4mm} R^{-1}_{1\to k}R^{-1}_{k\to 1}\\
&&+\sum_{s=2}^{k+1}\,M_{\,\ov 2}\dots \hat {M}_{\,\ov s}\dots M_{\,\ov{k+1}} \,R_1^{-1}\dots \RR_{s-1}\dots \RR_1,
\end{eqnarray*}
where the symbol $\hat{M}_{\,\ov s}$ means that  the term ${M}_{\,\ov s}$ is omitted. Note, that in the boundary term corresponding
to $s=2$ we assume $R^{-1}_1\dots R^{-1}_1\dots R^{-1}_1 = R^{-1}_1$.

As a consequence, we get the formula for the action of the operator $D_1$:
\be
D_1\triangleright M_{\,\ov 2}\dots M_{\,\ov{k+1}}=\sum_{s=2}^{k+1}\,M_{\,\ov 2}\dots \hat{M}_{\,\ov s}\dots M_{\,\ov{k+1}}\,R_1^{-1}\dots \RR_{s-1}\dots \RR_1,
\qquad \forall\, k\ge 1,
\label{act}
\ee
where we have taken into account that $D_1\triangleright 1=0$.
\begin{example} \rm A few first examples of (\ref{act}) are as follows:
$$
\begin{array}{l}
D_1\triangleright M_{\ov 2}= R_1^{-1}\\
\rule{0pt}{5mm}
D_1\triangleright M_{\ov 2}\, M_{\ov 3} = M_{\,\ov 3}\, R_1^{-1} + M_{\ov 2}\, R_1^{-1}R_2^{-1}R_1^{-1}\\
\rule{0pt}{5mm}
D_1\triangleright M_{\ov 2}\, M_{\ov 3}\, M_{\ov 4}=M_{\ov 3}\, M_{\ov 4}\,\RR_1+M_{\ov 2}\, M_{\ov 4}\,\RR_1\RR_2\RR_1+M_{\ov 2}\,
M_{\ov 3}\,\RR_1\RR_2\RR_3\RR_2\RR_1.
\end{array}
$$
 \end{example}

Now, we use the following formulas\footnote{In \cite{IOP} these formulas were established for a more general class of quantum
matrix algebras, associated with couples of
compatible braidings $(R,F)$. The  RE algebras correspond to the case $F=R$.}  \cite{IOP}:
\begin{eqnarray}
&& M_{\,\ov {k+1}}\dots M_{\,\ov {k+s}}=R_{k\,\to \,k+s-1}\, M_{\,\ov {k}}\dots M_{\,\ov {k+s-1}}\,R^{-1}_{k+s-1\,\to\, k}, \label{shift} \\
\rule{0pt}{5mm}&&R_iM_{\,\ov k} = M_{\,\ov k} R_i,\quad \forall\,i\not= k-1,k,\nonumber
\end{eqnarray}
which enable us to ``fill a gap'' $\hat M_{\,\ov s}$ in the right hand side of (\ref{act}) for any $2\le s\le k$:
$$
M_{\,\ov 2}\dots \hat{M}_{\,\ov s}\dots M_{\,\ov{k+1}}=R_{s\to k}\,M_{\,\ov 2}\dots M_{\,\ov {k}}\,\RR_{k\to s}.
$$
So, the action (\ref{act}) can be rewritten as:
\begin{eqnarray}
D_1\triangleright M_{\,\ov 2}\dots M_{\,\ov {k+1}}=\sum_{s=2}^{k} R_{s\to k}\, M_{\,\ov 2}
\dots M_{\,\ov {k}}&&\hspace*{-6mm} \RR_{k\to s}\,( \RR_1\dots \RR_{s-1}\dots \RR_1)
\label{act-nogap}\\
&&+M_{\,\ov 2}\dots M_{\,\ov {k}}\, \RR_1\dots \RR_{k}\dots \RR_1. \nonumber
\end{eqnarray}

Then, applying  formula (\ref{shift}) step by step $(m-1)$ times\footnote{\label{ftn:shift} Relation (\ref{shift}) is valid for the product 
of any matrix 
copies of the same size $A_{\,\ov k}B_{\,\ov{k+1}}\dots C_{\,\ov{k+p}}$, since it entirely follows from the braid relation for $R$ 
and does not depend on the algebraic properties of the 	matrices $A, B, \dots, C$.} we reduce (\ref{act-nogap}) to the following form:
\begin{eqnarray}
D_{\,\ov m}\triangleright M_{\,\ov {m+1}}\dots&&\hspace*{-6mm}M_{\,\ov {m+k}}=\nonumber\\
\sum_{s=m+1}^{m+k-1} &&\hspace*{-6mm}R_{s\to m+k-1}\, M_{\,\ov{m+1} }
\dots M_{\,\ov {m+k-1}}\,\RR_{m+k-1\to s}\, (\RR_m\dots \RR_{s-1}\dots \RR_m)\label{act-shift}\\
&&\hspace*{6mm}+M_{\,\ov {m+1}}\dots M_{\,\ov {m+k-1}}\, (\RR_m\dots \RR_{m+k-2}\dots \RR_m). \nonumber
\end{eqnarray}

At last, we make one more shift of the indices
$$
M_{\,\ov {m+1}}\dots M_{\,\ov {m+k-1}} = R_{m\rightarrow m+k-2}M_{\,\ov {m}}\dots M_{\,\ov {m+k-2}}\,R^{-1}_{m+k-2\rightarrow m}
$$
and present (\ref{act-shift}) in the form convenient for the subsequent application of $D_{\,\ov{m-1}}$:
$$
D_{\,\ov m}\triangleright M_{\,\ov {m+1}}\dots M_{\,\ov {m+k}} = \sum_{i=1}^kF_i(R)M_{\,\ov{m}}\dots M_{\,\ov {m+k-1}}\,G_i(R)
$$
where $F_i(R)$ and $G_i(R)$ are some chains of $R$-matrices $R_r$ with numbers $ r\geq m$ and consequently commuting with $D_{\,\ov p}$ if $p\leq m-1$.

Finally, the successive application of the operators $D_{\,\ov p}$ with $p=m-1,\dots ,1$ leads us to the following claim.

\begin{lemma} \label{lem:7} For any pair of positive integers $k\geq m$ the following relation takes place:
$$
D_{\,\ov{1\dots m}}\triangleright  M_{\,\ov {m+1}}\dots M_{\,\ov {m+k}}=\sum_j S_j(R) \, M_{\,\ov {1}}\dots M_{\,\ov {k-m}}\, T_j(R),$$
where $S_j(R)$ and $T_j(R)$ are some polynomials depending on $R_r,$ $1\le r\le m+k-1$.
If $k<m$, the result of this action vanishes.
\end{lemma}

Let us now fix an arbitrary element $z=Z(\tau_1,\dots ,\tau_{k-1})\in H_k(q)$ for some $k\ge 1$ and consider its image $ch_k(z)$ 
(see definition (\ref{sem})) in the characteristic
sublagebra of the RE algebra $\MM(R)$. The following matrix formula was proved in \cite{IOP}:
\be
I^{\otimes m}\,ch_k(z) = \Tr_{R(m+1\dots m+k)}\left(Z(R_{m+1},\dots ,R_{m+k-1})M_{\,\ov{m+1}}\dots M_{\,\ov{m+k}}\right).
\label{ch-sh}
\ee
With the use of this formula and relation in Lemma \ref{lem:7} we find the action of an arbitrary generalized Laplacian operator
$D_Q^{(m)}$ (\ref{gen-Lap}) on $ch_k(z)$:
 $$
D_Q^{(m)}\triangleright ch_k(z)=\Tr_{R(1\dots m+k)} \left(\sum_j Q(R) S_j(R)\,M_1\dots M_{\,\ov {k-m}}\, 
T_j(R)\, Z(R)\right),\quad k\ge m.
$$
Note that in virtue of the cyclic property of the quantum trace all polynomials in $R$ can be put together on the right (or left) hand side 
in this formula. 
Therefore, the result of action $D_Q^{(m)}\triangleright ch_k(z)$ belongs to the characteristic subalgebra of the RE algebra $\MM(R)$.
To complete the proof, we note that any element of the characteristic subalgebra of the RE algebra $\MM(R)$ is a finite sum of 
homogeneous
polynomials $ch_{k_i}(z_i)$ for some integers $k_i\ge 1$ and elements $z_i\in H_{k_i}(q)$.\hfill \rule{6.5pt}{6.5pt}

\section{Quantum Casimir operators}
\label{sec:4}

In this section we consider in detail the subalgebra $\hLL(R)$ of the quantum double $(\DD(\RR), \MM(R))$ generated by the linear
quantum differential operators
$$
\hl_i^j = \sum_{k=1}^N m_i^k\partial_k^j\,.
$$
The corresponding generating matrix is $\hL=MD$.

\begin{proposition}
The following claims hold true.
\begin{enumerate}
\item[\rm 1.]The matrix $\hL$ meets the quadratic-linear relation {\rm (\ref{mRE}):}
$$
R_1\hL_1 R_1\hL_1 - \hL_1R_1\hL_1 R_1 =  R_1\hL_1 -\hL_1 R_1.
$$
\item[\rm 2.] The permutation relations between generators $\hl_i^j$ and $m_k^s$ are given by the following permutation rules, 
presented in a matrix form as follows:\rm
\begin{equation}
R_1\hL_1R_1 M_1= M_1R_1 \hL_1  R^{-1}_1+R_1M_1.
\label{si}
\end{equation}\it
\end{enumerate}
\end{proposition}

\smallskip

\noindent
{\bf Proof. } Both claims of the proposition are straightforward consequences of the defining relations (\ref{RE}), (\ref{RE-D}) and
(\ref{raz}). Give a detailed proof of formula (\ref{si}).

Rewrite (\ref{raz}) in the equivalent form:
$$
D_1R_1M_1 = R_1M_1R_1^{-1}D_1R_1^{-1}+I^{\otimes 2}.
$$
Then relation (\ref{si}) can be obtained by the following chain of transformations:
$$
R_1\hL_1R_1 M_1 = R_1M_1\underline{D_1R_1 M_1} = \underline{R_1M_1R_1M_1}R_1^{-1}D_1R_1^{-1}+R_1M_1 = M_1R_1\hL_1R_1^{-1}+R_1M_1,
$$
where at the last step of transformations we used the relation (\ref{RE}) for $R_1M_1R_1M_1$.

The first claim can be proved in a similar way.\hfill \rule{6.5pt}{6.5pt}

\medskip

Introduce a set of central elements of the subalgebra $\hLL(R)$ which are homogeneous $k$-th order polynomials in generators
$\hl_i^j$:
\begin{equation} 
C^{(k)}_Q(\hL) = \Tr_{R(1\dots k)}\left(Q(R_1,\dots,R_{k-1})\, \hL_{\,\ov {1\dots k}}\right),\qquad \forall\,k\ge 1,
\label{gen-Cas}
\end{equation}
where $Q(R_1,\dots,R_{k-1})$ is an arbitrary polynomial in $R_i$, $1\leq i\leq k-1$.  Taking $Q = R_{k-1}\dots R_1$, we get the quantum
power sum $p_k(\hL)$ which can be rewritten in the form:
$$
p_k(\hL) =  \Tr_{R(1\dots k)}\left(R_{k-1}\dots R_1\, \hL_{\,\ov {1\dots k}} \right) = \Tr_R \hL^k.
$$

With respect to the action (\ref{D-M}) the central elements $C_Q^{(k)}(\hL)$ become the linear quantum differential operators
on the subalgebra $\MM(R)$. 

\begin{definition} \rm
The differential operators (\ref{gen-Cas}) will be called the generalized quantum Casimir operators. 
\end{definition}

The set of all
generalized Casimir operators $C^{(k)}_Q(\hL)$ for $\forall \,k\ge 1$ and arbitrary $Q(R)$ forms a central characteristic
subalgebra in $\hLL(R)$.

Our next aim is to investigate the action of the quantum differential operators $\hL$ on the RE algebra $\MM(R)$. In particular, we prove
that the characteristic subalgebra of $\MM(R)$ is mapped into itself by the action of any generalized Casimir operator $C^{(k)}_Q(\hL)$.

From a technical point of view it is more convenient to deal  with another set of $\hLL(R)$ generators $\hat K = \|\hat k_i^j\|$, the
corresponding generating matrices are connected by the following relation:
$$
\hat K = I - (q-q^{-1})\hL .
$$
It is easy to verify that  the matrix $\hat K$ satisfies the homogeneous quadratic relation (\ref{RE}):
$$
R_1\hat K_1 R_1\hat K_1 - \hat K_1R_1\hat K_1R_1 = 0.
$$

Let us fix the following notation:
\be
\hat K_{\underline 1} = \hat K_{1},\quad \hat K_{\,\underline {r+1}} = R^{-1}_r\hat K_{\,\underline r} R_{r}\quad \forall\,r\ge 1,\qquad \hat K_{\,\underline{r\rightarrow s}} =
\hat K_{\,\underline{r}}\hat K_{\,\underline{r-1}}\dots \hat K_{\,\underline {s+1}}\hat K_{\,\underline{s}}\quad r>s.
\label{K-cop}
\ee
In \cite{GPS6} it was proved that the action of $\hat K$ on an arbitrary monomial in generators $M$ of $\MM(R)$ can be written as
follows\footnote{See Proposition 10 of the cited paper.}:
\be
\hat K_{\,\underline{n+1}}\triangleright M_{\,\ov{1\rightarrow  n}} = J^{-1}_{n+1} M_{\,\ov{1\rightarrow  n}},
\label{K-act}
\ee
where $J_{n+1}^{-1} = R^{-1}_{n\rightarrow 1}R^{-1}_{1\rightarrow n}$ is the image of the inverse Jucys-Murphy element under the
$R$-matrix representation of the Hecke algebra $H_{n+1}(q)$. Note, that the action (\ref{K-act}) defines a representation of the RE
algebra $\hLL(R)$ in the algebra $\MM(R)$.

\begin{remark}\rm
In \cite{GPS6} the action (\ref{K-act}) was obtained in the frameworks of representation theory of the RE algebra. Alternatively, it can be
restored from (\ref{act-nogap}) by straightforward calculations.
\end{remark}

Taking into account that $\hat K_{\,\underline n}J_{m}^{-1} = J_{m}^{-1} \hat K_{\,\underline n}$ $\forall \, m<n$, we can subsequently
apply the formula (\ref{K-act}) to find the following general result:
\be
\hat K_{\,\underline{n+p\,\rightarrow\, n+1}}\triangleright M_{\,\ov{1\rightarrow  n}} =
\prod_{i=1}^{p}J^{-1}_{n+i}\prod_{s=2}^{p}J_s^{\,\uparrow \,n}\, M_{\,\ov{1\rightarrow  n}}.
\label{K-gen-act}
\ee
Here $J^{\,\uparrow \,n}_s = R_{n+s-1\,\rightarrow \, n+1}R_{n+1\,\rightarrow\, n+s-1}$ is the image of the Jucys-Murphy element
$j_s$ (\ref{JM}) under the $R$-matrix representation $\rho_R^{\,\uparrow\, n}$,  ``shifted'' by $n$ positions in the tensor product
of the spaces $V$. Thus, we have
\be
\rho_R^{\,\uparrow \,n}(\tau_i) = R_{n+i},\quad \Rightarrow \quad \rho_R^{\,\uparrow \,n}(j_s) = J_s^{\,\uparrow \,n} =
R_{n+s-1\,\rightarrow\, n+1}R_{n+1\,\rightarrow \,n+s-1}.
\label{rho-shif}
\ee

The following theorem is a direct consequence of the formula (\ref{K-gen-act}).
\begin{theorem}
\label{prop:10}
The action of any generalized quantum Casimir operator
$$
C^{(n)}_Q(\hat K) = \Tr_{R(1\dots n)}\left( Q(R_1,\dots, R_{n-1}) \,\hat K_{\,\underline{n\rightarrow 1}} \right)
$$
maps the characteristic subalgebra of the RE algebra $\MM(R)$ into itself.
\end{theorem}

\begin{remark}\rm
The form of $C^{(n)}_Q(\hat K) $ in the claim of the Theorem \ref{prop:10} fits well for our subsequent calculations. Actually it is
identical to the previous definition since for the
generating matrix of the RE algebra (\ref{RE}) the following identity holds true \cite{IP}:
$$
\hat K_{\,\underline {n\rightarrow 1}} = \hat K_{\,\ov{1 \rightarrow n}} \qquad \forall\,n\ge 1.
$$
This identity can be easily proved by induction in $n$.
\end{remark}

\smallskip

\noindent
{\bf Proof.} Choose an arbitrary homogeneous $n$-th order polynomial $ch_n(z)$ from the characteristic subalgebra of the RE algebra
$\MM(R)$
$$
ch_n(z) = \Tr_{R(1\dots n)}\left(Z(R_1,\dots R_{n-1})M_{\,\ov{1\rightarrow n}}\right)
$$
and prove that the action of an arbitrary generalized quantum Casimir operator $C^{(p)}_Q(\hat K)$ maps $ch_n(z)$ to an element of the characteristic subalgebra.
To do so we need the action (\ref{K-gen-act}) and the shift formula (\ref{ch-sh}) which we apply to the generalized Casimir operator:
$$
I^{\otimes n}\,C^{(p)}_Q(\hat K) = \Tr_{R(n+1\dots n+p)}\left(Q(R_{n+1},\dots ,R_{n+p-1})
\hat K_{\,\underline{n+p\,\rightarrow\, n+1}}\right).
$$
With the use of this formula, we present the action of the Casimir operator in the form:
\begin{eqnarray*}
C^{(p)}_Q(\hat K)\triangleright ch_n(z) &=& \Tr_{R(1\dots n)}\left(Z(R_1,\dots R_{n-1})C^{(p)}_Q(\hat K)\triangleright
M_{\,\ov{1\rightarrow n}}\right)\\
\rule{0pt}{5mm}&=& \Tr_{R(1\dots n+p)}\left(Z(R_1,\dots R_{n-1})Q(R_{n+1},\dots, R_{n+p-1})
\hat K_{\,\underline{n+p\,\rightarrow\, n+1}}\triangleright M_{\,\ov{1\rightarrow n}}\right).
\end{eqnarray*}
By  taking into account the action (\ref{K-gen-act}) we get the following answer:
\be
C^{(p)}_Q(\hat K)\triangleright ch_n(z) = \Tr_{R(1\dots n+p)}\left( F(R_1,\dots, R_{n+p-1}) M_{\,\ov{1\rightarrow n}}\right)
\label{C-M-pr}
\ee
where the polynomial $F$ reads:
$$
F(R_1,\dots,R_{n+p-1}) = Z(R_1,\dots R_{n-1})Q(R_{n+1},\dots, R_{n+p-1})\prod_{i=1}^{p}J^{-1}_{n+i}\prod_{s=2}^{p}J_s^{\,\uparrow\, n} .
$$
Using the permutation relations among the generators $\tau_i$ of the Hecke algebra $H_{n+p}(q)$ and the Hecke condition on $\tau_i$
one can show that any polynomial
$F(\tau_1,\dots,\tau_{n+p-1})$ can be presented in the form
\be
F(\tau_1,\dots,\tau_{n+p-1}) = Q_1(\tau_1,\dots, \tau_{n+p-2}) \tau_{n+p-1}Q_2(\tau_1,\dots, \tau_{n+p-2})+Q_3(\tau_1,\dots, \tau_{n+p-2})
\label{P-trans}
\ee
where the ploynomials $Q_i$, $i=1,2,3$, do {\it not} depend on $\tau_{n+p-1}$. So, in (\ref{C-M-pr}) we can consecutively transform the polynomials $F(R)$ under the $R$-trace as in (\ref{P-trans})
and calculate the $R$-traces in spaces with numbers form $n+p-1$ till $n+1$  using the property $\Tr_{R{(k+1)}}(R_{k}) = I_k$
(see (\ref{trR})). Thus, we come to the final result:
$$
C^{(p)}_Q(\hat K)\triangleright ch_n(z) = \Tr_{R(1\dots n)}\left(\tilde F(R_1,\dots, R_{n-1}) M_{\,\ov{1\rightarrow n}}\right).
$$
By definition, the right hand side lies in the characteristic subalgebra of the RE algebra $\MM(R)$ and, moreover, has the same
degree in $M$ as the initial polynomial $ch_n(z)$.

To complete the proof of the Theorem, we note that an arbitrary element of the characteristic subalgebra of $\MM(R)$ is a finite
sum of elements $ch_{n_i}(z_i)$ for some integers $n_i\ge 1$ and elements $z_i\in H_{n_i}(q)$. \hfill \rule{6.5pt}{6.5pt}

\section{Normal ordering and Wick theorem}
\label{sec:5}

In this section we introduce a quantum analog of normal ordering for products of  quantum first order differential operators $\hL = MD$.
We prove the Wick theorem for the product of partially ordered operators (Theorem \ref{th9}). In analogy with the classical case, 
the quantum normal ordering leads to operators, in which all quantum partial 
derivatives $\pa_i^j$ are placed to the right of ``coordinates'' $m_r^s$.

\begin{definition}\label{def:ord}\rm
The quantum normal ordering of derivatives $D = \|\partial_i^j\|$ and generators $M = \|m_i^j\|$ is defined by the rule:
\be
:\!D_1M_{\,\ov {2}}\!: \, =M_{\,\ov{2}}D_1 R_1^{-2}.
\label{razz}
\ee
\end{definition}
As usual, the ordered form of the product of any given elements $\partial_i^j$ and $m_k^r$ can be extracted form entries of matrix 
equality (\ref{razz}).

Note, that the rule (\ref{razz}) is valid in any double $(\DD(R^{-1}),\MM(R))$ defined with the use of a skew-invertible Hecke symmetry 
$R$, including the case of the supersymmetric $GL(m|n)$ type $R$-matrices. In the classical case $R=P$ where $P$ is the flip or 
super-flip the formula (\ref{razz}) reduces to the usual normal ordering of commutative
or supercommutative coordinates and corresponding partial derivatives.

Note that the transformation to the ordered form defined in (\ref{razz}) is performed with the permutation relations obtained from (\ref{raz}) by omitting the constant term:
\be
D_1 R_1M_1 = R_1M_1R_1^{-1}D_1R_1^{-1}.
\label{trans}
\ee
The general recipe of transformation of a product of quantum differential operators to the normal ordered form is as follows. Under the symbols of normal ordering
$:\,:$ one should apply the permutation relations (\ref{trans}) untill all quantum derivatives $D$ will be located on the right of all quantum ``coordinates'' $M$.

Let us give an example of such a transformation for the product of two quantum differential operators $\hL_1 = M_1D_1$ and 
$\hL_{\,\ov 2} = R_1\hL_1 R_1^{-1}$:
$$
:\!\hL_{1} \hL_{\,\ov 2}: \, = :\!M_1  \underline{D_1  R_1  M_1}  D_1 \RR_1\!:\, = M_1  R_1  M_1 \RR_1  D_1 \RR_1  D_1  \RR_1=
M_1\ M_{\,\ov 2} D_1  \RR_1  D_1 \RR_1.
$$
Taking into account that in the algebra $\DD(\RR)$  the following relation is valid:
$$
D_1  \RR_1  D_1=R_1  D_1  \RR_1  D_1  \RR_1= D_{\,\ov 2}  D_1  \RR_1,
$$
we get the final formula for the normal ordered form convenient for the subsequent generalization:
$$
:\!\hL_1 \hL_{\,\ov 2} :\,=M_1  M_{\,\ov 2}  D_{\,\ov 2}  D_1  R_1^{-2}=M_1  M_{\,\ov 2}  D_{\,\ov 2}  D_1 \JJ_2.
$$
With the use of permutation relations (\ref{raz}) we can express the product of two differential operators in terms of the ordered ones:
$$
\hL_1\, \hL_{\ov 2} =\, :\!\hL_1\, \hL_{\ov 2} : + \hL_1\, \RR_1.
$$

Our next aim is to determine the normal ordered form of an arbitrary order monomial in quantum operators $\hL$. To prove the corresponding theorem we need 
the following technical lemma.

\begin{lemma} 
\label{lem:order}
The normal ordered form of the operator $D_{\,\ov m}\hL_{\,\ov n}$ for $1\le m\le n-1$ is as follows:
\rm
\be
:\! D_{\,\ov m}\hL_{\,\ov n}:\, = \hL_{\,\ov n} D_{\,\ov m} \,J_{n-m}^{\,\uparrow m}(J_{n-m+1}^{-1})^{\uparrow (m-1)},\qquad 
n\ge 2,\quad 1\le m\le n-1,
\label{dl-norm}
\ee\it
where $J_{k}^{\,\uparrow p} = R_{k+p-1\,\rightarrow\, p+1}R_{p+1\,\rightarrow\, k+p-1}$ is the image of the Jucys-Murphy 
element $j_k$ under the shifted $R$-matrix 	repesentation $\rho_R^{\,\uparrow p}$ deined in {\rm (\ref{rho-shif})}.
\end{lemma}

\smallskip

\noindent
{\bf Proof.} For $m=1$, $n=2$ the claim of the lemma follows directly from (\ref{razz}):
$$
:\!D_1\hL_{\,\ov 2}:\,= :\!\underline{D_1M_{\,\ov 2}}D_{\,\ov 2}:\,= M_{\,\ov 2}D_1R_1^{-2}D_{\,\ov 2} = M_{\,\ov 2}D_1R_1^{-1}D_1R_1^{-1} =
M_{\,\ov 2}D_{\,\ov 2}D_1R_1^{-2} = \hL_{\,\ov 2}D_1 J_2^{-1}.
$$
Then, having applied an appropriate shift formula (\ref{shift}) $(m-1)$ times (see footnote \ref{ftn:shift}), we find
$$
:\! D_{\,\ov{m}}\hL_{\,\ov{m+1}}:\, =  \hL_{\,\ov{m+1}} D_{\,\ov{m}} R_{m}^{-2} = \hL_{\,\ov{m+1}} D_{\,\ov{m}} \,(J_{2}^{-1})^{\,\uparrow (m-1)},
$$
which coincides with (\ref{dl-norm}) for $n=m+1$ since  $J_1^{\,\uparrow m} = I_{m+1}$.

Now the final result (\ref{dl-norm}) for $n>m+1$ is easy to obtain:
\begin{eqnarray*}
:\! D_{\,\ov m}\hL_{\,\ov n}: &=& R_{n-1\rightarrow m+1}D_{\,\ov m}\hL_{\,\ov{m+1}}\,R^{-1}_{m+1\rightarrow n-1}
= R_{n-1\rightarrow m+1}\hL_{\,\ov{m+1}} D_{\,\ov m}\,R_m^{-2}R^{-1}_{m+1\rightarrow n-1}\\
\rule{0pt}{5mm}
&=& \hL_{\,\ov{n}} D_{\,\ov m}\,R_{n-1\rightarrow m+1}R_m^{-2}R^{-1}_{m+1\rightarrow n-1} = 
\hL_{\,\ov n} D_{\,\ov m} \,J_{n-m}^{\,\uparrow m}(J_{n-m+1}^{-1})^{\uparrow (m-1)}.
\end{eqnarray*}
The proof is completed.\hfill\rule{6.5pt}{6.5pt}

\medskip

We are able to prove the following main theorem on normal ordered forms.

\begin{theorem} 
\label{th7}
The normal ordered form of the product $\hL_{\,\ov{1\rightarrow k}} = \hL_1\dots \hL_{\,\ov k}$ reads as follows:\rm
\be
:\!\hL_{\,\ov{1\rightarrow k}}:\,=M_{\,\ov{1\rightarrow k}}D_{\,\ov{k\rightarrow 1}}\,
\Bigl(\prod_{s=1}^k\JJ_s\Bigr)=\Bigl(\prod_{s=1}^k\JJ_s\Bigr)M_{\,\ov{1\rightarrow k}}\,
D_{\,\ov{k \rightarrow 1}}, \qquad \forall\, k\ge1.
\label{fo}
\ee\it
\end{theorem}

\smallskip

\noindent
{\bf Proof.} We prove the theorem by induction in $k$. The base of induction at $k=1$ is obvious. Let (\ref{fo}) be valid up to some integer $k\ge 1$, we must verify
that then it is valid for $k+1$. 

Taking into account the induction assumption we write:
\be
:\!\hL_{\,\ov{1\rightarrow k+1}}: = \Big(\prod_{s=1}^k\JJ_s\Big):\!M_{\,\ov{1\rightarrow k}}\,D_{\,\ov{k \rightarrow 1}}\hL_{\,\ov{k+1}}:
\label{prom-res}
\ee
Then, the relation (\ref{dl-norm}) allows one to get the following normal ordered form:
$$
:\!D_{\,\ov {r\rightarrow 1}}L_{\,\ov{k+1}}: = L_{\,\ov{k+1}} D_{\,\ov {r\rightarrow 1}}\, J_{k-r+1}^{\,\uparrow r}J_{k+1}^{-1},\quad \forall\, r\le k.
$$
For $r=k$ this formula simplifies to $:\!D_{\,\ov {k\rightarrow 1}}L_{\,\ov{k+1}}: = L_{\,\ov{k+1}} D_{\,\ov {k\rightarrow 1}} J_{k+1}^{-1}$ 
and we can
complete our proof substituting this to (\ref{prom-res}):
$$
:\!\hL_{\,\ov{1\rightarrow k+1}}: =\Big( \prod_{s=1}^k\JJ_s \Big)M_{\,\ov{1\rightarrow k}}\,\hL_{\,\ov{k+1}}\,D_{\,\ov{k \rightarrow 1}}
J_{k+1}^{-1} =
\Big( \prod_{s=1}^{k+1}\JJ_s \Big)M_{\,\ov{1\rightarrow k+1}}\,D_{\,\ov{k+1 \rightarrow 1}}.
$$
Here we used the fact that matrix $J^{-1}_i$ commute with the product $M_{\,\ov{1\rightarrow p}}\,D_{\,\ov{p\rightarrow 1}}$ for all 
$i\le p$ due to relations
$$ 
R_s^{-1}  M_{\,\ov s}   M_{\,\ov{s+1}}=M_{\,\ov s}  M_{\,\ov{s+1}} R_s^{-1}, \qquad R_s^{-1}  D_{\,\ov{s+1}} D_{\,\ov{s}}=
D_{\,\ov{s+1}}   D_{\,\ov{s}}  R_s^{-1}.
$$
Also, this allows one to place the multipliers $\JJ_s$ at any side in formula (\ref{fo}).\hfill\rule{6.5pt}{6.5pt}
\medskip

Now we are going to establish a quantum analog of the Wick formula, which transforms a partially ordered product of differential 
operators into the totally ordered form. To prove the quantum Wick theorem we need the lemma below.

\begin{lemma}
\label{lem:dl}
The following matrix identity takes place for $\forall\, k\ge 1$:\rm
\be
D_{\,\ov{ k\,\rightarrow\,1}}\,\hL_{\,\ov{k+1}}=\hL_{\,\ov{k+1}}\, D_{\,\ov{ k\,\rightarrow \,1}}J_{k+1}^{-1}+
D_{\,\ov{ k\,\rightarrow \,1}}\,\frac{I_{k+1}-J^{-1}_{k+1}}{q-\qq}.
\label{induc}
\ee \it
\end{lemma}

\smallskip

\noindent
{\bf Proof.} We prove the lemma by induction in $k$. The induction base for $k=1$ is an immediate consequence of (\ref{raz}) and (\ref{RE-D}). Indeed, multiplying
(\ref{raz}) by $R_1^{-1}D_1R_1$ from the right and applying (\ref{RE-D}) for $R_1^{-1}D_1R_1^{-1}D_1$ we obtain:
\be
D_1\hL_{\,\ov 2} = \hL_{\,\ov 2}D_1 R_1^{-2} + D_1R_1^{-1},
\label{D-L}
\ee
which is precisely formula (\ref{induc}) written for $k=1$ if we take into account the definition of the Jucys-Murphy element:
$$
J_2^{-1} = R_1^{-2} = I_2-(q-q^{-1})R_1^{-1}.
$$
On applying relation (\ref{shift}) successively $k$ times to (\ref{D-L}) we extend it to the higher matrix copies:
\be
D_{\,\ov k}\hL_{\,\ov{k+1 }} = \hL_{\,\ov{k+1}}D_{\,\ov k} R_k^{-2} + D_{\,\ov k}\,R_k^{-1}.
\label{D-L-h}
\ee

Now, assume that formula (\ref{induc}) is valid up to some integer $k-1\ge 1$. We should prove that then it is fulfilled for $k$ too.
We have a chain of transformations:
\begin{eqnarray*}
D_{\,\ov{ k\,\rightarrow\,1}}\,\hL_{\ov{k+1}} &= &D_{\,\ov{ k\,\rightarrow\,1}}\,R_k\,\hL_{\ov{k}}\,\RR_k=D_{\ov k}\,R_k\,
\un{D_{\,\ov{ k-1\,\rightarrow\,1}}\,\hL_{\ov{k}}}\,\RR_k\\
\rule{0pt}{5mm} &=&
D_{\ov k}\,R_k\,\hL_{\ov k}\,D_{\,\ov{ k-1\,\rightarrow\,1}} \,\JJ_k\,\RR_k+ D_{\,\ov{ k\,\rightarrow\,1}} \, R_k\, 
\frac{I_k-\JJ_{k}}{q-\qq}\,\RR_k\\
\rule{0pt}{5mm} &=&  D_{\ov k}\,\hL_{\ov{k+1}}\,D_{\,\ov{ k-1\,\rightarrow\,1}} \,R_k\, \JJ_k\,\RR_k+ 
D_{\,\ov{ k\,\rightarrow\,1}} \,  \frac{I_{k+1}-R_k\,\JJ_{k}\,\RR_k}{q-\qq}.
\end{eqnarray*}
Here, the underlined term is transformed in accordance with the induction assumption.

Next, we use (\ref{D-L-h}) and continue the  above transformations as follows:
\begin{eqnarray*}
D_{\,\ov{ k\,\rightarrow\,1}}\,\hL_{\ov{k+1}} &= & \hL_{\,\ov{k+1}}\,  D_{\,\ov{ k\,\rightarrow\,1}}\,\RR_k\, \JJ_k\, \RR_k+ 
D_{\,\ov{ k\,\rightarrow\,1}} \left( \JJ_k\,R_k^{-1}+
\frac{I_{k+1}-R_k\, \JJ_k\, \RR_k}{q-\qq}\right)\\
\rule{0pt}{5mm} &=&
\hL_{\,\ov{k+1}}\, D_{\,\ov{ k\,\rightarrow\,1}}\,\JJ_{k+1}+D_{\,\ov{ k\,\rightarrow\,1}}\,\frac{I_{k+1}-\JJ_{k+1}}{q-\qq}.
\end{eqnarray*}
The proof of the lemma is completed.\hfill\rule{6.5pt}{6.5pt}

\medskip

So, we are ready to prove a quantum analog of the Wick theorem.

\begin{theorem} \label{th9}
The following quantum Wick formula holds true:\rm
\be
:\!\hL_{\,\ov{1\rightarrow k}}:\hL_{ \,\ov{k+1}}=\,:\!\hL_{\,\ov{1\,\rightarrow\, k+1}}:+:\!\hL_{\,\ov{1\rightarrow k}}:
\frac{I_{k+1}-J^{-1}_{k+1}}{q-\qq},\qquad \forall\,k\ge 1.
\label{the9}
\ee\it
\end{theorem}

\smallskip

\noindent
{\bf Proof.}
In accordance with Theorem \ref{th7} we have:
$$
:\!\hL_{\,\ov{1\rightarrow k}}:\, \hL_{\ov{k+1}}=\Big(\prod_{s=1}^k\JJ_s\Big) M_{\,\ov{1\,\rightarrow\, k}}\, D_{\,\ov{ k\,\rightarrow 1}}\,\hL_{\ov{k+1}}.
$$
Taking into account (\ref{induc}) we get:
\begin{eqnarray*}
\Big(\prod_{s=1}^k\JJ_s\Big) M_{\,\ov{1\,\rightarrow\, k}}\,D_{\,\ov{k\,\rightarrow\, 1}}
\hL_{\ov{k+1}}&\hspace*{-2mm}=&\hspace*{-3mm}\Big(\prod_{s=1}^k\JJ_s\Big)
M_{\,\ov{1\,\rightarrow\, k}}\, \hL_{\ov{k+1}}\,D_{\,\ov{k\,\rightarrow\, 1}} \,\JJ_{{k+1}}\\
\rule{0pt}{5mm}&\hspace*{-5mm}+&\hspace*{-4mm}\Big(\prod_{s=1}^k\JJ_s\Big) M_{\,\ov{ 1\,\rightarrow\,k}}\, 
D_{\,\ov{ k\,\rightarrow 1}}\,\frac{I_{k+1}-\JJ_{k+1}}{q-\qq}\\
\rule{0pt}{5mm}&\hspace*{-5mm}=&\hspace*{-4mm}\Big(\prod_{s=1}^k \JJ_s\Big)M_{\,\ov{ 1\,\rightarrow\,k+1}}\, 
D_{\,\ov{ k+1\,\rightarrow 1}}\,\JJ_{{k+1}} +
:\!\hL_{\,\ov{1\,\rightarrow\,k}}:\,\frac{I_{k+1}-\JJ_{k+1}}{q-\qq}\\
\rule{0pt}{5mm}&\hspace*{-5mm}=&\hspace*{-4mm}:\!\hL_{\,\ov{1\,\rightarrow\,k+1}}:+
:\!\hL_{\,\ov{1\,\rightarrow\,k}}:\,\frac{I_{k+1}-\JJ_{k+1}}{q-\qq}.
\end{eqnarray*}
The proof is completed.\hfill\rule{6.5pt}{6.5pt}

\section{Universal quantum Capelli identities}
\label{sec:6}

In this section we establish the set of {\it universal quantum matrix Capelli  identities}. The term ``universal'' reflects the fact, that
all other known forms of Capelly identities (see \cite{Ok1, Ok2, JLM}) can be obtained as particular cases (as matrix projections or 
as a limit $q\to 1$) of the universal ones.

Let us introduce the following matrix notation:
\begin{equation}
\PP_1=I_1,\qquad \PP_k=\frac{I_{k}-J^{-1}_k}{q-\qq}, \quad k\geq 2.
\label{P-mat}
\end{equation}
Taking into account the definition of $J_k^{-1}$ and the Hecke condition $R_{\,i}^{-2} = I_{i+1}-(q-q^{-1})R_{\,i}^{-1}$ we present 
$\PP_{k+1}$ as the following polynomial in $R_{\,i}$, $ i\leq k$:
$$
\PP_{k+1}=\RR_k+\sum_{s=1}^{k-1}\RR_{k\,\rightarrow\,s+1}\,R_s^{-1} \RR_{s+1\,\rightarrow k}.
$$
This form of $\PP_{k+1}$ is convenient for calculating the classical limit $q\rightarrow 1$. Observe that if the Hecke symmetry $R$ 
tends to the flip $R\rightarrow P$ as 
$q\rightarrow 1$, the polynomial $\PP_{k+1}$ tends to the sum of the transpositions $(i,k+1)$ in the tensor representation of the 
group algebra ${\Bbb C}[S_{k+1}]$ in the space $V^{\otimes (k+1)}$:
$$
\PP_{k+1}\,\rightarrow \,\sum_{i=1}^{k}P_{ik+1}.
$$

\begin{theorem}
For the quantum differential operators $\hL = MD$ the following matrix Capelli identities take place for $\forall\, k\ge 2$:\rm
\be
\hL_1(\hL_{\ov 2}-\PP_2)\dots (\hL_{\ov k}-\PP_k)=M_{\,\ov{1\,\rightarrow \,k}}\, D_{\,\ov{ k\,\rightarrow \,1}}\,\Big(\prod_{s=1}^k \JJ_s\Big).
\label{MatCap}
\ee
\end{theorem}

\smallskip

\noindent
{\bf Proof.} To prove (\ref{MatCap}) it suffices to rewrite (\ref{the9}) in the form:
$$
:\!\hL_{\,\ov{1\,\rightarrow \,k-1}}:\left(\hL_{\ov{k}}-\PP_{k}\right) = :\!\hL_{\,\ov{1\,\rightarrow\,k}}:
$$
and then apply the same formula for $:\!\hL_{\,\ov{1\,\rightarrow \,k-1}}:$ and so on. For $:\!\hL_{\,\ov{1\,\rightarrow\,k}}:$ in the right 
hand side we use (\ref{fo}).
\hfill\rule{6.5pt}{6.5pt}

\medskip

If the Hecke symmetry $R$ is a deformation of the usual flip $P$ or of the superflip of $GL(m|n)$ type, then passing to the limit $q\rightarrow 1$, we get the corresponding universal 
matrix Capelli identities in $U(gl(N))$ or in $U(gl(m|n))$.

Now we point out an important corollary of the Capelli identities (\ref{MatCap}): the normal ordering (\ref{razz}) maps the characteristic
subalgebra of the RE algebra $\hLL(R)$ into itself. In other words, the normal ordering does not destroy
the centrality of a central differential operator. 
\begin{theorem}\label{theor:22}
The normally ordered form of any generalized Casimir operator {\rm (\ref{gen-Cas})}
$$
:\!C^{(k)}_Q(\hL)\!: \, = \Tr_{R(1\dots k)}\left(Q(R_1,\dots, R_k):\!\hL_{\,\overline{1\to k}}\!:\right)
$$
belongs to the characteristic subalgebra of the RE algebra $\hLL(R)$ and, therefore, is a central quantum differential operator.
\end{theorem}

\smallskip

\noindent
{\bf Proof.} As follows from (\ref{fo}) and (\ref{MatCap}), the normally ordered form of $\hL_{\,\overline{1\to k}}$ coincides with
the left hand side of the corresponding Capelly identity:
$$
:\!\hL_{\,\overline{1\to k}}\!: \,= \hL_1(\hL_{\ov 2}-\PP_2)\dots (\hL_{\ov k}-\PP_k).
$$
Then, taking into account the definition (\ref{P-mat}) it is not difficult to verify that
$$
\PP_s \hL_{\,\overline r} = \hL_{\,\overline r}\,\PP_s, \quad \forall\, r>s.
$$
So, having expanded all the brackets in the right hand side of the above expression for $:\!\hL_{\,\overline{1\to k}}\!:$, we can
move all the matrices $\PP_s$ to the {\it right} of all chains of $\hL$ operators and get the following expression:
\begin{equation}
:\!\hL_{\,\overline{1\to k}}\!: \, = \hL_{\,\overline{1\to k}} + \sum_{n=1}^{k-1}\,\,\sum_{2\le s_1<\dots<s_n\le k}
\hL_1\dots \hL'_{s_1}\dots \hL'_{s_n}\dots \hL_{\,\overline{k}}\,\PP_{s_1}\dots \PP_{s_n},
\label{L-exp}
\end{equation}
where the symbols with prime mean the absence of the corresponding multipliers.

Next, we use the formula (\ref{shift}) in order to ``fill the gaps" in the product of operators $\hL$,  and present the typical
term in (\ref{L-exp}) in the form (see the proof of the Theorem \ref{theor:8}):
$$
\hL_1\dots \hL'_{s_1}\dots \hL'_{s_n}\dots L_{\,\overline{k}} = F_{\{s_1,\dots, s_n\}}(R)\,\hL_{\,\overline{1\,\to \,k-n}}\,
F_{\{s_1,\dots, s_n\}}(R^{-1})
$$
for some polynomials $F_{\{s_1,\dots, s_n\}}(R)$.  

And, at last, the cyclic property of the $R$-trace allows us to come to the final answer:
\begin{equation}
:\!C^{(k)}_Q(\hL)\!:  = C^{(k)}_Q(\hL) + \sum_{n=1}^{k-1}C^{(k-n)}_{Q^{(n)}}(\hL),
\label{ord-Cas}
\end{equation}
where the polynomials $Q^{(n)}(R)$ are found as partial $R$-traces of the form:
$$
Q^{(n)}(R) = \sum_{2\le s_1<\dots<s_n\le k}\,\Tr_{R(k-n+1\dots k)}\left(
F_{\{s_1,\dots, s_n\}}(R^{-1})\,\PP_{s_1}\dots \PP_{s_n}Q(R)F_{\{s_1,\dots, s_n\}}(R)\right).
$$
So, any normally ordered generalized Casimir operator can be presented as a finite sum of generalized Casimir operators (\ref{ord-Cas}) 
and, therefore, it is a central quantum differential operator from the characteristic subalgebra.\hfill\rule{6.5pt}{6.5pt}

\medskip

In conclusion, we would like to explain the connection of the universal Capelli identities with the results, obtained in \cite{JLM}, 
\cite{Ok1} and \cite{Ok2}.
If we multiply (\ref{MatCap}) from any side by the image $E^\lambda_T(R) = \rho_R(e^{\lambda}_T)$ of the primitive idempotent
 $e^{\lambda}_T$, $\lambda\vdash k$ of the Hecke 
algebra $H_k(q)$ we find the identity obtained in \cite{JLM}, theorem 4.1. This fact is a direct consequence of the formula (see, 
for example, \cite{OP}):
$$
J_s^{-1}E^\lambda_T(R) = E^\lambda_T(R)J_s^{-1} = q^{-2c_s(T)} E^\lambda_T(R),\quad \forall\,\lambda\vdash k,\quad1\le s\le k,
$$
where $c_s(T)$ is the content of the box with an integer $s$ in the standard Young table $T$ of the Young diagram $\lambda$. 
Then, on taking the limit $q\rightarrow 1$ for $R$ being 
the deformation of the usual flip, we come to the identities in $U(gl(N))$ obtained in \cite{Ok1,Ok2}.

Emphasize, that our version of the Capelli identities (\ref{MatCap}) does not depend on projectors 
$E^\lambda_T(R)$ and it is valid for any skew-invertible Hecke symmetry $R$, including the supersymmetric $GL(m|n)$ 
type $R$-matrices.

\end{document}